\documentclass{article}
\usepackage[latin1]{inputenc}
\usepackage[small]{titlesec}
\usepackage{amsmath}
\usepackage{amsthm}
\usepackage{amssymb}
\usepackage{graphicx}
\usepackage{standalone}
\usepackage{color}
\theoremstyle{plain}
\newtheorem{thm}{Theorem}[section]

\newtheorem{rem}[thm]{Remark}

\usepackage[small]{titlesec}
\usepackage{blindtext}

	\title{List of constructions of $NO^+(6,2)$}
\author{Valentino Smaldore \footnote{Valentino Smaldore:
valentino.smaldore@unipd.it 
Dipartimento di Tecmica e Gestione dei Sistemi Industriali,
Universit\`{a} degli Studi di Padova, Stradella San Nicola 2, 36100 Vicenza, Italy.}}

\date{}

\begin{document}
\maketitle
\begin{abstract}
In this note, we list many equivalent constructions of the tangent polar graph $NO^+(6,2)$.

\end{abstract}

A \textit{strongly regular graph} with parameters $(v,k,\lambda,\mu)$ is a graph with $v$ vertices where each vertex is incident with $k$ edges, any two adjacent vertices have $\lambda$ common neighbours, and any two non-adjacent vertices have $\mu$ common neighbours. We denote it by $srg(v,k,\lambda,\mu)$. It is easy to show that the complement of the strongly regular graph $srg(v,k,\lambda,\mu)$ is a $srg(v,v-k-1,v-2k+\mu-2, v-2k+\lambda)$. We provide above a list of constructions, equivalent up to isomorphisms, for a strongly regular graph with $v=28$, $k=15$, $\lambda=6$, $\mu=10$.

\section{$NO^+(6,2)$}\label{sec1}
Let $Q^+(2n-1,2)$ be a non-degenerate hyperbolic quadric of $PG(2n-1,2)$. Let $NO^+(2n, 2)$ be the graph whose vertices are the points of $PG(2n-1,2) \setminus Q^+(2n-1,2)$ and two vertices $P_1$ and $P_2$ are adjacent if the line joining $P_1$ and $P_2$ is a line tangent to $Q^+(2n-1, 2)$. The graph $NO^+(2n, 2)$ is a strongly regular graph with parameters $(2^{2n-1}-2^{n-1},2^{2n-2}-1,2^{2n-3}-2,2^{2n-3}+2^{n-2})$, see \cite{brvan} for more detailed information. When $n=3$, the \textit{Klein Quadric} $Q^+(5,2)$ is the set of points of the variety given by:
\begin{equation*}
    Q: X_1X_6+X_2X_5+X_3X_4=0.
\end{equation*}
The quadric contains 35 points, and the generators of the quadric are the 30 subplanes, which split into two families of size 15, called respectively \textit{Latin planes} and \textit{Greek planes}. The graph $NO^+(6,2)$ is a $srg(28,15,6,10)$. The automorphism group is isomorphic to the orthogonal group $PGO^+(6,2)$, stabilizing the Klein Quadric in $PGL(6,2)$.

\section{Quadric with a hole}\label{sec2}
An alternative description was provided in \cite{BIK} and was later analyzed in a geometrical setting by \cite{CDMPS}. Consider a hyperbolic quadric $Q^+(2n-1,q)$ in $PG(2n-1,q)$, and let $\Pi$ be a generator of the quadric, i.e. a maximal totally isotropic subspace. Hence, we define the graph $\overline{\mathcal{G}_n}$ by taking as vertices the points in $X=Q^+(2n-1,q)\setminus\Pi$, adjacent if the projective line joining them is contained in $X$. The strongly regular graph obtained has parameters $v=q^{n-1}\left(\frac{q^n-1}{q-1}\right)$, $k=q^{n-1}\left(\frac{q^{n-1}-1}{q-1}\right)$, $\lambda=q^{n-1}\left(\frac{q^{n-2}-1}{q-1}\right)+q^{n-2}(q-1)$, $\mu=q^{n-1}\left(\frac{q^{n-2}-1}{q-1}\right)$. If $q=2$ the graph $\overline{\mathcal G_3}$ is a $srg(28,12,6,4)$, while its complement $\mathcal G_3$ is isomorphic to $NO^+(6,2)$, see \cite{RS} for the explicit isomorphism. The graph $\mathcal G_n$ can be described as the graph on $X$, with two adjacency relations: let $P,Q \in X$, then $P \sim_1 Q$ if and only if the line $\langle P,Q \rangle$ is secant to $Q^+(2n-1,q)$, and
$P \sim_2 Q$ if and only if the line $\langle P,Q \rangle$ is contained in $Q^+(2n-1,q)$ and meets the generator $\Pi$ in a point. 
\begin{rem}
    Note that while $q=2$, the graph $\mathcal G_n$ always has the same set of parameters as $NO^+(2n,2)$, but the two graphs are isomorphic if and only if $n\leq3$, see \cite{CDMPS}. While $n\geq4$, the automorphism group of $\mathcal G_n$ is isomorphic to $Stab_{PGO^+(2n,2)}(\Pi)$.
\end{rem}

\section{Klein correspondence} \label{sec3}
We use $L$ to denote the set of all lines of $PG(3,q)$. For any $\ell=\langle x,y\rangle\in L$, with $x=(x_{0}, x_{1}, x_{2}, x_{3})$, $y=(y_{0}, y_{1}, y_{2}, y_{3})$, set
$$p_{ij}=\begin{vmatrix} x_i & x_j\\
y_i & y_j \end{vmatrix},$$
$i, j=0, 1, 2, 3$. The map
$$\mathcal{K}:\,
   \begin{cases}
   L\rightarrow PG(5,q)\\
   \langle x,y\rangle\mapsto \langle(p_{12},p_{13},p_{14},p_{23},p_{24},p_{34})\rangle.
   \end{cases}$$
is  called \textit{Klein correspondence}, and it maps the lines in $L$ to points in the Klein Quadric $Q^+(5,q)$, see \cite{2}.

 \begin{table}[h!]
 \centering
 \footnotesize{\begin{tabular}{|c|c|}
\hline
\textbf{Lines of $PG(3,q)$} & \textbf{Points of $Q^+(5,q)$}\\
\hline
Two skew lines & Two non-orthogonal points\\
\hline
Two intersecting lines & Two orthogonal points\\
\hline
$q+1$ lines on a plane-pencil  &  $q+1$ points on a line\\
\hline
$q^{2}+q+1$ lines through a point & Latin plane $\mathcal{L}$\\
\hline
$q^{2}+q+1$ lines on a plane & Greek plane $\mathcal{G}$\\
\hline
\end{tabular}}
\caption{\footnotesize{Images of structures of $PG(3,q)$ under the action of the Klein map $\mathcal{K}$.}}
\end{table}
Now, if the generator $\Pi$ is a Greek plane, we describe the graph as follows. Consider a projective space $PG(3,2)$ and fix a Fano subplane $\pi=PG(2,2)$. Hence, the vertex set of the graph will be made of the lines in $L$ not lying on $\pi$, by fixing two adjacency rules: $\ell$ and $r$ are adjacent if they either do not intersect ($\sim_1$) or intersect in $\pi$ ($\sim_2$). Example given, let $\ell\cap\pi=\{P\}$, then $\ell$ has 12 neighbours in $\sim_1$ as there are $35-3\cdot6-1=16$ lines not meeting $\ell$ in $PG(3,2)$, 4 of them in $\pi\setminus\{P\}$; and $\ell$ has 3 neighbours in $\sim_2$, since such neighbours are the lines not in $\pi$ meeting $\ell$ in $P$.\\ In an analogous way, if the generator $\Pi$ of $Q^+(5,2)$ is a Latin plane, we consider a projective space $PG(3,2)$ and fix a point $Q$. The vertex set of the graph will consist of the lines in $L$ that do not meet $Q$, by fixing two adjacency rules: $\ell$ and $r$ are adjacent if they either do not intersect ($\sim_1$) or lie in $\langle\ell,Q\rangle$ ($\sim_2$).

\section{Antiflags in $PG(2,2)$}\label{sec4}
Recently, in \cite{PI} another alternative construction for the graphs $NO^+(2n,2)$ has been shown, taking as vertices the point-hyperplane antiflags in $PG(n-1,2)$, i.e. a point $P$ and a hyperplane $\Pi$ such that $P\notin\Pi$. This provides a 4-class association scheme for antiflags $(P,\Pi)$ with the following relations.
\begin{itemize}
    \item[$A_0$:] $(P,\Pi)\sim_0(P',\Pi')$ if $P=P'$ and $\Pi=\Pi'$;
    \item[$A_1$:] $(P,\Pi)\sim_1(P',\Pi')$ if $P\in\Pi'$ and $P'\notin\Pi$ or $P\notin\Pi'$ and $P'\in\Pi$;
    \item[$A_2$:] $(P,\Pi)\sim_2(P',\Pi')$ if $P\in\Pi'$ and $P'\in\Pi$;
    \item[$A_3$:] $(P,\Pi)\sim_3(P',\Pi')$ if either $P=P'$ and $\Pi\neq\Pi'$ or $P\neq P'$ and $\Pi=\Pi'$;
    \item[$A_4$:] $(P,\Pi)\sim_4(P',\Pi')$ if $P\notin\Pi'$ and $P'\notin\Pi$.
\end{itemize}
We obtain a graph isomorphic to $NO^+(2n,2)$ taking the relations $A_2,A_3,A_4$. In particular, while $n=3$ the graph $NO^+(6,2)$ has vertex set consisting of the point-line antiflags $(P,\ell)$ in the projective plane $PG(2,2)$, adjacent if in relation $A_i$, $i=2,3,4$.\\
Since non-degenerate conics in $PG(2,2)$ are complements of point-line antiflags, we provide an easy alternative description considering as vertices of $NO^+(6,2)$ the 28 non-degenerate conics of the projective plane over the binary field, while considering the same adjacency rules. The latter construction will have a deeper interpretation considering the Veronese embedding of conics in $PG(2,q)$ into $PG(5,q)$.

\section{$N\mathcal M_4^3$}\label{sec5}
The \textit{Veronese surface} of all conics of $PG(2,q)$, is the variety $\mathcal{V}^{4}_{2}=\{(a^{2},b^{2},c^{2},ab,ac,bc)|(a,b,c)\in PG(2,q)\}\subseteq PG(5,q)$. The mapping
  $$\mu:\,
   \begin{cases}
    PG(2,q)\rightarrow PG(5,q)\\
   (x_{1},x_{2},x_{3})\mapsto(x_{1}^{2},x_{2}^{2},x_{3}^{2},x_{1}x_{2},x_{1}x_{3},x_{2}x_{3}).
   \end{cases}$$
   is called the Veronese embedding of $PG(2,q)$. The variety $\mathcal{V}^{4}_{2}$ consists of $q^{2}+q+1$ points and its stabilizing subgroup of $\mathcal{V}^{4}_{2}$ in $PGL(6,q)$ arises by \textit{lifting} from the collineation group of $PG(2,q)$.  The group of lifted collineations has the following orbits on the $q^{5}+q^{4}+q^{3}+q^{2}+q+1$ conics of $PG(2,q)$:
   \begin{itemize}
    \item $\mathcal{O}_{1}:=q^{2}+q+1$ double lines (points of $\mathcal{V}^{4}_{2}$);
    \item $\mathcal{O}_{2}:=\frac{1}{2}(q^{2}+q+1)(q^{2}-q)$ pairs of imaginary lines;
    \item $\mathcal{O}_{3}:=\frac{1}{2}(q^{2}+q+1)(q^{2}+q)$ pairs of intersecting lines;
    \item $\mathcal{O}_{4}:=q^{5}-q^{2}$ non-degenerate conics.
   \end{itemize}
   The set of all degenerate conics $\mathcal{O}_{1}\cup\mathcal{O}_{2}\cup\mathcal{O}_{3}=\mathcal{M}^{3}_{4}$ is called \textit{secant variety} and  $|\mathcal{M}^{3}_{4}|=|Q^{+}(5,q)|=(q^{2}+1)(q^{2}+q+1)$, see \cite[Theorem 4.18]{2}. The secant variety $\mathcal{M}^{3}_{4}$ is a hypersurface of degree 3 and dimension 4. We may identify points of $PG(5,q)$ with $3\times3$ symmetric matrices in $S_3(\mathbb F_q)$, by:
   $$(X_{1},X_{2},X_{3},X_{4},X_{5},X_{6})\longleftrightarrow\left(\begin{array}{ccc}
                                                                                        X_{1} & X_{4} & X_{5} \\
                                                                                        X_{4} & X_{2} & X_{6} \\
                                                                                        X_{5} & X_{6} & X_{3}
                                                                                      \end{array}\right).$$
   In this representation, the Veronese surface $\mathcal{V}^{4}_{2}$ correspond to the matrices $\left(\begin{array}{ccc}
                                                                                        x_{1}^{2} &x_{1}x_{2} & x_{1}x_{3} \\
                                                                                        x_{1}x_{2} & x_{2}^{2} & x_{2}x_{3} \\
                                                                                        x_{1}x_{3} & x_{2}x_{3} & x_{3}^{2}
                                                                                      \end{array}\right),$ while $\mathcal{M}^{3}_{4}$ is a cubic hypersurface with equation
   \begin{equation}
    \left|\begin{array}{ccc}
     X_{1} & X_{4} & X_{5} \\
     X_{4} & X_{2} & X_{6} \\
     X_{5} & X_{6} & X_{3}
    \end{array}\right|=0.
   \end{equation}
   With the above notation, the orbit $\mathcal{O}_{1}=\mathcal{V}^{4}_{2}$ coincides with the $3\times3$ symmetric matrices over $GF(q)$ of rank 1, while $\mathcal{O}_{2}$ and $\mathcal{O}_{3}$ with the $3\times3$ symmetric matrices over $GF(q)$ of rank 2, and $\mathcal{O}_{4}$ with the $3\times3$ symmetric matrices over $GF(q)$ of rank 3. When $q=2$, the set $PG(5,2)\setminus\mathcal{M}^{3}_{4}$ has cardinality 28, as $|V(NO^{+}(6,2))|$, because $|Q^{+}(5,2)|=|\mathcal{M}^{3}_{4}|=35$. 
   Since $\mathcal{M}^{3}_{4}\cap\mathcal{K}=\mathcal{O}_{2}=N$, the secant variety always shares a plane with a Klein Quadric. Hence, we can give an alternative construction of the graph $NO^+(6,2)$, which here is called $N\mathcal{M}^{3}_{4}$ as in \cite{RS}.
    \begin{itemize}
     \item $V(\mathcal{M}^{3}_{4})=PG(5,2)\setminus\mathcal{M}^{3}_{4}$;
     \item $E(\mathcal{M}^{3}_{4})=\{(x,y)|x,y\in V(\mathcal{M}^{3}_{4}),|\langle x,y\rangle\cap \mathcal{M}^{3}_{4}|=1\}.$
    \end{itemize}

\section{Nonsingular $3\times3$ matrices}\label{sec6}
Finally, it is also possible to describe the graph in the representation as $3\times3$ matrices over $\mathbb{F}_{q}$:
   \begin{itemize}
     \item $V(\mathcal{M}^{3}_{4})$ is the set of the non-singular symmetric matrices of order 3 over $\mathbb{F}_{q}$;
     \item $E(\mathcal{M}^{3}_{4})=\{(A,B)|A,B\in V(\mathcal{M}^{3}_{4}),A+B \mbox{ is singular }\}.$
    \end{itemize}

    \section{Conclusion}
    We end up summarizing all the constructions provided, focusing on the vertex set of the graph and its ambient space.
    
 \begin{table}[h!]
 \centering
 \footnotesize{\begin{tabular}{|c|c|c|c|}
\hline
 & \textbf{Graph} & \textbf{Ambient space} & \textbf{Vertex set}\\
\hline
Section \ref{sec1} & $NO^+(6,2)$ & $PG(5,2)$ & $PG(5,2)\setminus Q^+(5,2)$\\
\hline
Section \ref{sec2} & $\mathcal G_3$ & $PG(5,2)$ & $Q^+(5,2)\setminus\Pi$, $\Pi$ generator of the quadric\\
\hline
Section \ref{sec3} & $NO^+(6,2)$ & $PG(3,2)$ &  Lines in $PG(3,2)$ not on a fixed plane\\
\hline
Section \ref{sec3} & $NO^+(6,2)$ & $PG(3,2)$ &  Lines in $PG(3,2)$ not through a fixed point\\
\hline
Section \ref{sec4} & $NO^+(6,2)$ & $PG(2,2)$ & Point-line antiflags $(P,\ell)$, $P\notin\ell$.\\
\hline
Section \ref{sec4} & $NO^+(6,2)$ & $PG(2,2)$ & Non-degenerate conics\\
\hline
Section \ref{sec5} & $N\mathcal M^3_4$ & $PG(5,2)$ & $PG(5,2)\setminus \mathcal M_4^3$\\
\hline
Section \ref{sec6} & $N\mathcal M^3_4$ & $S_3(\mathbb 
Z_2)$ & $M\in S_3(\mathbb Z_2)$ such that $det(M)\neq0$\\
\hline
\end{tabular}}
\end{table}

\end{document}